\documentclass[10pt]{article}
\usepackage[tbtags]{amsmath}
\usepackage{amssymb,color}
\usepackage{cite}

\def\polhk#1{\setbox0=\hbox{#1}{\ooalign{\hidewidth
    \lower1.5ex\hbox{`}\hidewidth\crcr\unhbox0}}}
\definecolor{c20}{rgb}{0.,0.7,0.}
\definecolor{c30}{rgb}{0.,0.,1.}
\definecolor{c40}{rgb}{1,0.1,0.7}
\definecolor{c50}{rgb}{1,0,0}

\def\pE#1{\textcolor{c40}{#1}}
\def\pE#1{#1}
\def\hH#1{\textcolor{c40}{#1}}
\def\hH#1{#1}
\def\hh#1{\textcolor{c30}{#1}}
\def\hh#1{#1}

\def\rH#1{\textcolor{c40}{#1}}
\def\rH#1{#1}
\def\aH#1{\textcolor{c40}{#1}}
\def\aH#1{#1}
\def\bH#1{\textcolor{c40}{#1}}
\def\bH#1{#1}

\def\cL#1{\textcolor{c30}{#1}}
\def\cL#1{#1}

\def\cc#1{\textcolor{c50}{#1}}
\def\cc#1{#1}

\def\lcx#1{\textcolor{c50}{#1}}
\def\lcx#1{#1}

\catcode`\@=11 \@addtoreset{equation}{section} \catcode`\@=12

\newcommand{\nwc}{\newcommand}
\nwc{\COM}[1]{}

\newtheorem{theo}{Theorem}[section]
\newtheorem{sat}[theo]{Proposition}
\newtheorem{de}[theo]{Definition}
\newtheorem{lem}[theo]{Lemma}
\newtheorem{korr}[theo]{Corollary}
\newtheorem{remark}[theo]{Remark}
\newcommand{\nelem}[1]{{Lemma \ref{#1}}}

\newcommand{\netheo}[1]{{Theorem \ref{#1}}}

\newtheorem{exxa}[theo]{Example}

\newcommand{\abs}[1]{\lvert #1 \rvert}
\newcommand{\Abs}[1]{ \Bigl \lvert #1 \Bigr \rvert}

\newcommand{\pk}[1]{{\mathbb{P}} \{#1\} }

\newcommand{\Pk}[1]{\mathbb{P}\left \{#1 \right\}}
\newcommand{\R}{\mathbb{R}}

\newcommand{\inr}{\in \R}

\newcommand{\limit}[1]{\lim_{#1 \to   \infty}}

\newcommand{\todis}{\ \stackrel{d}{\to} \ }
\newcommand{\toprob}{ \stackrel{p}{\to}}

\newcommand{\equaldis}{\stackrel{d}{=}}

\newcommand{\BAN}{\begin{align}}
\newcommand{\EAN}{\end{align}}
\newcommand{\BANY}{\begin{align*}}
\newcommand{\EANY}{\end{align*}}
\newcommand{\BQN}{\begin{eqnarray}}
\newcommand{\EQN}{\end{eqnarray}}
\newcommand{\BQNY}{\begin{eqnarray*}}
\newcommand{\EQNY}{\end{eqnarray*}}
\newcommand{\BS}{\begin{sat}}
\newcommand{\ES}{\end{sat}}
\newcommand{\BL}{\begin{lem}}
\newcommand{\EL}{\end{lem}}
\newcommand{\BT}{\begin{theo}}
\newcommand{\ET}{\end{theo}}
\newcommand{\BK}{\begin{korr}}
\newcommand{\EK}{\end{korr}}
\newcommand{\BD}{\begin{de}}
\newcommand{\ED}{\end{de}}
\newcommand{\BIT}{\begin{itemize}}
\newcommand{\EIT}{\end{itemize}}
\newcommand{\BDI}{\begin{description}}
\newcommand{\EDI}{\end{description}}

\newcommand{\QED}{\hfill $\Box$}
\newcommand{\IF}{\infty}

\def\fracl#1#2{\left( \frac{#1}{#2} \right) }

\newcommand{\prooftheo}[1]{ \textsc{Proof of Theorem} \ref{#1}:}

\newcommand{\prooflem}[1]{\textsc{Proof of Lemma} \ref{#1}:}

\begin{document}
%%%%%%%%%%%%%%%%%%%%%%%%%%%%%%%%%%%%%%%%%%%%%%%%%%%%%%%%%

\centerline{\Large \bf Maxima of Skew Elliptical Triangular \hH{Arrays}}

       \vskip 0.8 cm
       \centerline{\large Enkelejd Hashorva and Chengxiu Ling}

       \vskip 0.2 cm

      \centerline{University of Lausanne, 1015 Lausanne, Switzerland}

      \centerline{\today{}}

       \vskip 1 cm

{\bf Abstract:} In this paper we investigate the asymptotic
behaviour of the componentwise maxima for two  bivariate skew
elliptical triangular arrays with components given in terms of skew
\pE{transformations} of bivariate spherical random vectors. We find
the weak limit of the normalized maxima \pE{for both cases that} the
random radius pertaining to the elliptical random vectors is either
in the Gumbel or in the Weibull max-domain \pE{of} attractions.

{\it Key words and phrases:} \pE{maxima} of triangular arrays; skew
elliptical distributions; Gumbel max-domain attraction; Weibull
max-domain attraction; H\"usler-Reiss distribution;
max-infinitely divisible distribution.

\section{Introduction}\label{sec1}
In the seminal paper \cite{HuslerR89}, {H\"{u}sler and Reiss showed that }the maxima of dependent Gaussian triangular arrays \pE{converge} \hH{in distribution} (after normalization)
to a random vector with max-stable distribution function \pE{(df)} referred to as H\"usler-Reiss \pE{df}. Specifically,
if $(X_{jn},Z_{jn}), j\le n, n\ge 1$  is a triangular array of independent bivariate Gaussian random vectors
such that $(X_{jn}, Z_{jn}) \equaldis (X, \rho_n X +
\sqrt{1-\rho_n^2} Y)$ for $\rho_n \in (-1,1)$ where $X,Y$ are independent $N(0,1)$ random variables with \pE{df} $\Phi$
and  $\equaldis$ stands for equality of \pE{dfs},
then the convergence in distribution
\BQN \label{eqA}
 \Biggl(  \frac{\max_{1 \le j \le n} X_{jn}- b_n}{a_n}, \frac{\max_{1 \le j \le n} Z_{jn}- b_n}{a_n}\Biggr)
\todis (\mathcal{M}_1, \mathcal{M}_2), \quad n\to \IF \EQN holds
with $a_n= 1/\sqrt{2 \ln n}, b_n=\Phi^{-1}(1-1/n) $ and
$(\mathcal{M}_1, \mathcal{M}_2)$ having the H\"usler-Reiss
\pE{df} $H_\lambda$ given by \BQNY \label{H} H_\lambda
(x, y) = \exp \left( - \Phi \left( \lambda
+\frac{x-y}{2\lambda}\right) e^{-y} - \Phi \left( \lambda
+\frac{y-x}{2\lambda}\right) e^{-x}\right), \quad \lambda \in
\lcx{[0,\IF)},\EQNY provided that the H\"usler-Reiss
condition
\BQN\label{HRcond}
 \limit{n}(1- \rho_n)  \ln n &=& \lambda^2
 \EQN holds. \pE{The marginal dfs of $H_\lambda$ are the unit Gumbel df $\Lambda(x)= \exp(-\exp(-x)), x\inr$.}

\pE{A popular extension of the normal distribution is the skew normal one, see e.g., \cite{Azzalini2003, Ma2004}};
a natural question \hH{is} \cc{whether} skewing has an effect on the asymptotic behaviour of
the normalized maxima. \pE{Specifically, define}
a triangular array $(X_{jn},Z_{jn}), j\le n, n\ge 1$  of \pE{bivariate skew normal random vectors as }
 \BQN
\cc{(X_{jn}, Z_{jn}) \equaldis (\abs{X}, \rho_n \abs{X} +
\sqrt{1-\rho_n^2} Y), \quad \rho_n \in (-1,1), \quad j\le n.} \EQN
If the H\"usler-Reiss condition \eqref{HRcond} holds, then by conditioning it follows that %\pE{(see Appendix)}
 \eqref{eqA} \pE{is still valid} where instead of
\cc{ $b_n = \Phi^{-1}(1-1/n)$, we {put} $b_n = \Phi^{-1}(1-1/n) +
a_n\ln2$}; in view of the findings of \cite{ChiDeHLK,HashorvaKW12}
the above result is expected. However, for the half-skew
model, surprisingly the \pE{limiting} \pE{df}  is \cc{not \hH{the}
H\"usler-Reiss \pE{df} any more}. \pE{Indeed,}  if we suppose
that for all $j\le n$
 \BQN \label{H-skew}
 (X_{jn}, Z_{jn}) \equaldis (X, \rho_n \abs{X} +
\sqrt{1-\rho_n^2} \cc{Y}), \quad \rho_n \in (-1,1), \EQN then again
under the H\"usler-Reiss \rH{condition} \pE{\eqref{HRcond}} we obtain the joint convergence
in distribution \BQN \label{eqAB}
 \Biggl(  \frac{\max_{1 \le j \le n} X_{jn}- b_n}{a_n},
 \cc{\frac{\max_{1 \le j \le n} Z_{jn}- (b_n + a_n \ln 2)}{a_n}}
 % \frac{\max_{1 \le j \le n} Z_{jn}- b_n}{a_n+ \ln 2}
 \Biggr)
\todis (\widetilde{\mathcal{M}_1}, \widetilde{\mathcal{M}_2}), \quad
n\to \IF\cc{,}
\EQN
where $(\widetilde{\mathcal{M}_1},
\widetilde{\mathcal{M}_2})$ has \pE{df}
$\widetilde{H_\lambda}$ given by
\BQN
 \widetilde{H_\lambda}(x,y)=  H_{\lambda}(x, y+ \ln 2) \Lambda (y+ \ln 2), \quad x,y\inr.
 \EQN
As shown for instance in \cite{HB} the H\"usler-Reiss \pE{df} \rH{appears due to}
the fact that Gaussian random vectors belong to the larger class of
elliptically symmetrical  {(elliptical for shorthand)} %(for short elliptical)
random vectors with pertaining random radius which has \pE{df}
in the Gumbel max-domain \pE{of} attraction (MDA), see definition below.\\
In this paper we are concerned with the behaviour of maxima \aH{of the larger
classes of} skew and half-skew elliptically symmetric distributions.
\pE{For simplicity, we shall deal with the bivariate} setup assuming that
$(X,Y)$ has stochastic representation
\BQN\label{Def.1}
 (X, Y) \equaldis R (\cos \Theta, \sin \Theta),
 \EQN
where the random \lcx{angle} $\Theta$ is independent of $R>0$ and
follows \rH{the uniform} distribution \cL{on} $(-\pi, \pi)$
\rH{(abbreviate this by $\Theta \sim U(-\pi, \pi)$)}. In the
\rH{special} case that $R^2$ is chi-square distributed with 2
degrees of freedom then $X$ and $Y$ are
independent \pE{with df} \cc{$\Phi$}. % $F$.

Our investigation will be concerned with two different cases:  {$a)$} %i)
we consider %the cases that are similar to the Gaussian one mentioned
{the case that the random radius $R$ has df in
the Gumbel \pE{MDA} which includes the Gaussian one
mentioned above,}
and  b) %ii)
 we shall assume
that $R$ has df in the Weibull \pE{MDA}, see definition below. The motivation for the latter
assumption comes from \cite{Hashorva08,HB}. Interestingly, as shown
in the aforementioned papers the case of \pE{the} Weibull \pE{MDA} leads to a limiting distribution \hH{which is different
from the} H\"usler-Reiss distribution, see also \cite{HB}.  This will be confirmed
in this paper also for the skew elliptical and half-skew elliptical models.

\COM{ So far there are no papers in the literature considering
extreme value behavior of skew elliptical triangular array.
Therefore, the aim of this paper is to establish the limit laws of
the skew elliptical triangular array. Specifically, for given
$\rho_n, \rho_{1n}, \rho_{2n} \in (0, 1]$, define two skew
elliptical triangular arrays of independent bivariate random vectors
$(X_{jn}, Z_{jn})$ and $( Z_{jn}^{(1)}, Z_{jn}^{(2)}), j = 1,
\ldots, n$ via the stochastic representation (cf. \cite{AzzaliniD96,
Fang03})

\BQN\label{Def.2} (X_{jn}, Z_{jn}) \equaldis (X, \rho_n |X| +
\sqrt{1-\rho_n^2} Y), \quad (Z_{jn}^{(1)}, Z_{jn}^{(2)}) \equaldis
(\rho_{1n} |X| + \sqrt{1-\rho_{1n}^2} Y, \rho_{2n} |X| +
\sqrt{1-\rho_{2n}^2} Y), \EQN

with $(X, Y)$ defined in \eqref{Def.1}. If $\rho_n, \rho_{1n},
\rho_{2n}$ goes to 1 as $n\to\IF$, then ($\toprob$ denotes
convergence in probability)

\[(X_{jn}, Z_{jn}) \toprob (X, |X|), \quad (Z_{jn}^{(1)},
Z_{jn}^{(2)}) \toprob (|X|, |X|), \quad \forall j\ge1. \]

This may eventually imply asymptotic dependence of the sample maxima
given by
\[(M_{n1}, M_{n2}) = (\max_{1\le j\le n} X_{jn},
\max_{1\le j\le n} Z_{jn}), \quad (\widetilde {M_{n1}}, \widetilde
{M_{n2}}) = \Big(\max_{1\le j\le n} Z^{(1)}_{jn}, \max_{1\le j\le n}
Z^{(2)}_{jn}\Big). \]

In the present paper we shall investigate the asymptotic behavior of
the triangular array defined in \cc{\eqref{H-skew}} assuming that
the associated random radius $R$ is both in the Gumbel and in the
Weibull \pE{MDA}, respectively. }

Our main findings are presented in \netheo{T1} and \netheo{T2} below. Therein we show that maxima of triangular arrays of
half-skew elliptically symmetric  distributions is different as {those} % that
of elliptically symmetric and skew elliptically symmetric distributions. %Specifically,
%we shall prove that both the norming constants and
%the limiting distributions are  different, which is not the case for the skew elliptically distributions. \\
\pE{Since for the skew elliptical case the limiting distributions are new our findings are of some theoretical and applied interest;
in a forthcoming paper we shall discuss aspects of statistical modelling using skew and half-skew elliptically symmetric  distributions.}
%The flexibility offered by the large class of skew and half-skew elliptical distributions in various statistical modelling One one hand
%However, the norming constants are not the same. Furthermore, we show that for the case of half-skew elliptical distributions both the limiting distributions and the norming constants are different that The main
%contributions of this paper \rH{concern the} derivation of the
%normalizing constants such that the maxima of the skew \rH{and half-skew} elliptical
%triangular \rH{arrays is} in the \pE{MDA}. \rH{Furthermore, we show} that the asymptotic dependence of the sample maxima may
%eventually depend on the convergence of \cc{$\rho_n \to 1$ as in
%\cite{HashorvaKW12, HuslerR89}} according to the different skew extent.

% \rH{Organization of the paper:}
\pE{The organisation of the paper is as follows:} \aH{In Section
\ref{sec2} we present our main results.} \rH{Two} illustrating
examples are given in \cc{Section} \ref{sec3}. The proofs of the results are relegated to Section \ref{sec4}.
%{We conclude the paper with a short Appendix.}

\section{Main Results} \label{sec2}
 In this section % for $(X, Y)$ specified by  \eqref{Def.1}
 we shall investigate the asymptotic behavior of maxima of triangular arrays
 defined via %\aH{(1.2) and (1.3)}
 \eqref{H-skew} % and  {\eqref{Def.3} below}
 with $(X, Y)$ specified by  \eqref{Def.1}.
 The main {assumption} imposed on the random radius $R$ is that it has
df $F$ in the Gumbel  or in the Weibull \pE{MDA}.
%In view of Lemma 2.1 in \cite{HashorvaPT10}
There is no loss of generality if we fix the right endpoint $x_F$ of $F$. Therefore, hereafter
we shall assume that $x_F \in \{{1},\IF\}$. % \{0,\IF\}$.

\subsection{Gumbel MDA}
Next we deal with the case that $F$ is in the Gumbel \pE{MDA} \hh{with some positive scaling function $w$,} i.e.,
 \hH{for any $s\inr$} %\hh{(set $\bar F = 1-F$)}
\BQN \label{Gumbel} \lim_{u\uparrow  x_F} \frac{\bar F(u + s/w(u))}{\bar
F(u)} = \exp(-s), \quad \bar F: = 1-F. %x \in \R.  %\hH{\bar F= 1- F},
 \EQN
The scaling function $w$ can be defined
asymptotically by \lcx{(cf. \cite{Embrechts97, HRF})}
\[w(t) =\pE{(1+o(1))}  \frac{\bar F(t)}{\int_t^{x_F} \bar F(s)\, ds}, \quad t\to x_F.\]
For $X= R \cos\Theta$ as defined in \eqref{Def.1} the assumption
\eqref{Gumbel} implies that the df  $G$ of $X$ is also in the Gumbel
MDA with the same scaling function $w$ as {for} $F$. In fact also
the converse result holds, see Theorem 4.1 in \cite{MR2678878}.
Furthermore,
from the aforementioned theorem or Theorem 12.3.1 in \cite{Berman92} (see also Theorem 3 in \cite{Polar92}, Proposition A.3 in \cite{convex} %
%Theorem 4.1 in \cite{regina},
 and Theorem 1 in \cite{HashExt12} for more general results)
%$B(\cdot; a, b)$ be the df of a Beta random variable with positive
%parameters $a,b$.
%\BL \label{L0}
%(Theorem 12.3.1 in \cite{Berman92})
%Suppose $F$ belongs to the Gumbel \pE{MDA} with scaling
%function $w$ satisfying
condition \eqref{Gumbel} implies
with $v(u) = \sqrt{u
w(u)}$
\BQNY \int_0^1 \bar F\fracl{u}{\sqrt{1-y}}\, dB(y; a,
b)=
 \frac{2^a\Gamma(a+b)}{\Gamma(b)}(v(u))^{-2a} \bar F(u)
 (1+o(1)), \quad u\uparrow x_F,
 \EQNY
where  $B(y;a,b)$ is the beta df with positive parameters $a,b$ which has density
$$ \frac{\Gamma(a+b)}{\Gamma(a)\Gamma(b)} x^{a-1}(1- x)^{b-1}, \quad x\in (0,1).
$$
Here  {$\Gamma(\cdot)$ denotes the Euler Gamma function.
%{\remark \label{R0}  Let $X$ be given by \eqref{Def.1} with df $G$.
Since $(\cos\Theta)^2$ has beta df $B(y;1/2, 1/2)$ we have  thus
%distribution with parameteIn our particular case, Then in view of \nelem{L1} with $a=b=1/2$, $G$ is in the Gumbel
%\pE{MDA} \hH{with the same scaling function $w$. Specifically,
\BQN \label{F and G}
\nonumber \aH{2 \pk{X> u}}= \Pk{ R|\cos \Theta|
> u}=\sqrt{ \frac 2 \pi} \frac{ \bar F(u) }{\nu(u)} (1 + o(1)), \quad u\uparrow x_F.
%\quad \mathrm{with}.\ \nu(u) := \sqrt {u w(u)}.
\EQN
This fact can be also formulated in the classical framework \pE{of extreme value theory} as
\[ \lim_{n\to\IF} \sup_{x\inr} \Abs{ G^n(a_n x + b_n) - \Lambda(x)}=0,
\]
 with
 \BQN \label{bn}  b_n  =
G^{\leftarrow} (1-1/n) = \inf\{x: G(x)
> 1-1/n\}, \quad a_n = \frac 1{w(b_n)}.
 \EQN

Next we state our first result.
 \BT \label{T1}  {Let $(X, Y)$ be given by \eqref{Def.1} with $F$ the df of the random radius $R$ \pE{which satisfies \eqref{Gumbel}
 with some positive scaling function $w$.}  Let $\rho_n \in (0, 1]$ and let $a_n, b_n$ be given by
 \eqref{bn}.}
%  Suppose that $F$ is in the Gumbel \pE{MDA} satisfying \eqref{Gumbel} with auxiliary function
% $w$ and
If %Suppose further that
\BQN\label{Cond.1} \limit{n}\frac{(1 - \rho_n) b_n }{a_n} =
2\lambda^2 \in \lcx{[0, \IF)}\EQN
and \pE{ $(X_{jn},Z_{jn}), j\le n, n\ge1$ is given by
 \eqref{H-skew}}, then %$a)$  For $Z_{jn}, j\le n, n\ge1$
%given by \eqref{H-skew} }
\BQN\label{T1r1} \frac{\max_{1 \le j \le n} Z_{jn}- (b_n + a_n \ln
2)}{a_n} \todis \mathcal{M}, \quad n\to \IF, \EQN where
$\mathcal{M}$ has {unit Gumbel} df $\Lambda$. Furthermore, the
\aH{joint convergence in \eqref{eqAB} holds}.
 \ET

{\remark \label{R1} $a)$ \aH{In the case that  $\rho_n  \pE{=}\rho \in (0,
1)$ for all large $n$}, then it follows  that \netheo{T1} holds with
$\lambda = \IF$ and $\widetilde {H_\IF}(x, y) =\Lambda (x)
\Lambda (y).$  \\
$b)$ Let $(\widetilde {M_{n1}}, \widetilde {M_{n2}}) = (\max_{1 \le
j \le n} Z_{jn}^{(1)}, \max_{1 \le j \le n} Z_{jn}^{(2)} )$ with
% $(Z_{jn}^{(1)}, Z_{jn}^{(2)}) \equaldis (\rho_{1n} |X| +
% \sqrt{1-\rho_{1n}^2} Y, \rho_{2n} |X| + \sqrt{1-\rho_{2n}^2} Y)$
\BQN \label{Def.3}
 (Z_{jn}^{(1)}, Z_{jn}^{(2)}) \equaldis (\rho_{1n} |X| +
\sqrt{1-\rho_{1n}^2} Y, \rho_{2n} |X| + \sqrt{1-\rho_{2n}^2} Y).
\EQN Suppose that $\rho_{\cL{i}n}$ satisfies condition
\eqref{Cond.1} with $\lambda_i \in [0,\IF), \cL{i} =1, 2$. \pE{Using
\nelem{L1} and following the arguments of the proof of  \netheo{T1}
we obtain}
 \BQNY \left( \frac{\widetilde {M_{n1}}- (b_n + a_n \ln
2)}{a_n}, \frac{\widetilde {M_{n2}}- (b_n + a_n \ln 2)}{a_n} \right)
\todis (\mathcal M_1, \mathcal M_2),  \quad n\to \IF,% \quad \mathrm{\lcx{with}} \quad
% \lambda = |\lambda_1 -  \lambda_2|
 \EQNY
{where $(\mathcal M_1, \mathcal M_2)$ has H\"usler-Reiss
df $H_\lambda$ with $\lambda = |\lambda_1 -
 \lambda_2|$}%as in the Introduction
 , see \eqref{H}.\\
c) \bH{The bivariate H\"usler-Reiss df appeared initially
in \cite{bro1977}; see also the recent articles {\cite{eng2011,
% kab2011,
kab2009, oes2012}}. Related results for more general
triangular arrays can be found in {\cite{eng2012b,DavisF,
frick2010limiting, HashorvaW13, kab2011,manjunath2012some}}; see also
\cite{eng2012a} for novel statistical applications}}.

\subsection{Weibull case}
The main assumption on $F$ in this section is that it belongs to the Weibull MDA, hence $x_F$ is necessarily finite. We assume
again that $x_F=1$. Specifically, we shall \aH{suppose} that $F$ is
in the Weibull \pE{MDA}  with index $\alpha>0$, i.e., (cf. \cite{Embrechts97})
\BQN\label{Weibull}
\lim_{t \downarrow 0}\frac{\bar F( {1} -tx)}{\bar F(1-t)} = x^{\alpha},
\quad x>0.
\EQN
In view of Theorem 4.5 in \cite{MR2678878} %also
the df $G$ of $X$ is in the Weibull MDA, and the converse result is also valid.
For notational simplicity  \pE{define below a positive constant $\mathcal{I}_\alpha$} and a df $\Upsilon_\alpha (\cdot)$ }
by
 \BQN \label{I(alpha)}
 \mathcal {I}_\alpha = \frac{\sqrt2}{\pi} \int_{0}^1(1-s^2)^\alpha \, ds = \frac 1{\sqrt{2\pi}} \frac{\Gamma(\alpha +1)}{\Gamma(\alpha + 3/2)},
  \quad \Upsilon_\alpha (t) = \frac{\int_{-1}^ t (1-s^2)^\alpha \, ds}{\int_{-1}^ 1 (1-s^2)^\alpha \,
 ds},\quad t \in [-1, 1].
 \EQN
\BT \label{T2} Let $(X,Y)$ be given by \eqref{Def.1} with $F$ the df of the random radius $R$ and
$\rho_n\in(0, 1]$. Suppose that $F$ {is in the Weibull \pE{MDA} satisfying \eqref{Weibull}} for
some $\alpha
> 0$.

$a)$ If for some $ u_n \downarrow 0$ \BQN \label{Cond.11}
\limit{n}\frac{1 - \rho_n }{ u_n } = 2\lambda^2 \in [0, \IF), \EQN
then
\BQN\label{Res.00} \limit{n} \frac{\Pk{\rho_n |X| + \sqrt{1 -
\rho_n^2} Y
> 1 + u_n x } } { u_n
^{1/2} \bar F(1 - u_n) } = 2 \mathcal {I}_\alpha  \abs x^{1/2 +
\alpha} , \quad x<0.
 \EQN
$b)$ {Let $(X_{jn}, Z_{jn}), j\le n, n\ge1$ given by \eqref{H-skew}.
If \eqref{Cond.11} holds for $u_n$ such that $u_n ^{1/2} \bar F(1 -
u_n) =\pE{(1+o(1))} /n$, then}
 \BQN \label{Limit.1}& & \limit{n}\Pk{ \frac{\max_{1 \le j \le n} X_{jn}- 1}{\bH{u_n}} \le x,
 \frac{\max_{1 \le j \le n} Z_{jn}- \bH{1}}{\bH{u_n}} \le y } \notag\\
& \quad & = \exp\left( - \mathcal {I}_\alpha \left( \abs x ^{1/2 +
\alpha} \Upsilon_\alpha \fracl{\lambda + (y-x) /
(2\lambda)}{\sqrt{2\abs x}} + \abs y^{1/2 + \alpha}\left( 1 +
\Upsilon_\alpha \fracl{\lambda + (x-y) / (2\lambda)}{\sqrt{2\abs y}}
\right) \right) \right)
 \EQN
 {holds for all $x, y <0$.}% and $a_n = u_n, b_n = 1$.}
\ET

{\remark \label{R2}  {a)} The right-hand side of
\eqref{Limit.1} is a bivariate df with dependent
Weibull
marginal \pE{dfs}.\\
$b)$ If $\rho_n \pE{=}\rho \in (0, 1)$ for all large $n$, then $\lambda =
\IF$ and the triangular array $(X_{jn}, Z_{jn}), j \le n, n\ge1$  given by \eqref{H-skew} is
asymptotically tail
independent. \\
$c)$ {Let $(\widetilde {M_{n1}}, \widetilde {M_{n2}})$ be given as in Remark \ref{R1} $b)$. % With the same notation as in Remark \ref{R1} $b)$ we obtain
Suppose that $\rho_{in}$ satisfies condition \eqref{Cond.11} with
$\lambda_i \in [0,\IF), i =1, 2$ and $u_n$ such that $u_n^{1/2} \bar
F (1-u_n) =\pE{(1+o(1))} /n$.} Using \nelem{L1} and following the
arguments of the proof of  \netheo{T2} we obtain that for all $x, y
< 0$
 \BQNY
 & & \limit{n}\Pk{ \frac{\widetilde {M_{n1}}- 1}{\bH{u_n}} \le x, \frac{\widetilde {M_{n2}}- 1}{\bH{u_n}} \le
 y }  \\
& \quad & = \exp\left( -2 \mathcal {I}_\alpha \left( \abs x ^{1/2 +
\alpha} \Upsilon_\alpha \fracl{\lambda + (y-x) /
(2\lambda)}{\sqrt{2\abs x}} + \abs y^{1/2 + \alpha} \Upsilon_\alpha
\fracl{\lambda + (x-y) / (2\lambda)}{\sqrt{2\abs y}} \right)
\right)
 \EQNY
holds with $\lambda = \abs{\lambda_1 - \lambda_2}$. {Note that the
limit joint df above is the same as that without skew
transformation, see (2.8) therein in Theorem 2.1 in
\cite{Hashorva08}. }
% with $\lambda = \abs{\lambda_1 - \lambda_2}$ {provided that there exists a positive sequence $u_n$ such that} $u_n^{1/2}
% \bar F (1-u_n) =\pE{(1+o(1))} /n$ and
% \BQNY \limit{n}\frac{1 - \rho_{jn} }{u_n} =
% 2\lambda_j^2 \in [0, \IF), \quad j = 1, 2.\EQNY

$d)$ For both cases that $F$ is  {in} the Gumbel or in the Weibull
\pE{MDA} \pE{utilising the arguments in \cite{Hashorva08} the more
general cases that $(\sin \Theta)^2$ has  regularly varying tail at $\cL{\pi/2}$ %1
can be dealt with following the ideas of our proofs
here; we omit those results.} }
\section{Examples}\label{sec3}

In this section, we give  {two} illustrating examples.

{\exxa \label{Exam.1} (Kotz Type $I$)} Consider a triangular array
as in \netheo{T1} with almost surely positive random radius $R$ \bH{and df  $F$ such that as $x\to \IF$}
\[\bar F(x) =\pE{(1+o(1))} K x^{\varsigma} \exp\left( - c x^\tau \right), \quad \varsigma \in \R, K, c, \tau >0.\]
\bH{It follows} that $F$ is in the Gumbel \pE{MDA} with
auxiliary function $w(x) =c \tau x^{\tau -1}$ and \lcx{further} by
\eqref{F and G}
\[ \bar G(x) {= (1+o(1))}\frac{K}{\sqrt{2\pi c \tau}} x^{\varsigma - \tau/2} \exp\left( - c x^\tau \right) =:  (1+o(1))\widetilde K x^{\widetilde \varsigma } \exp\left( - c x^\tau \right), \quad x\to \IF,\]
which implies that  {(see e.g., \cite{Embrechts97})}% \aH{(see e.g. Embrechts et al.\ (1997))}
\[b_n = G^\leftarrow (1-1/n) %\sim (c^{-1}\ln n)^{1/\tau} + \frac1{c\tau}(c^{-1}\ln n)^{1/\tau -1} \left(\frac{\widetilde \varsigma }{\tau} \ln(c^{-1} \ln n) + \ln \widetilde K  \right)
 { = (1+o(1))} (c^{-1}\ln n)^{1/\tau}, \quad n\to \IF \]
 and thus $b_n / a_n  = (1+o(1)) c \tau ( b_n )^\tau  = (1+o(1)) \tau \ln n$.
 Therefore, condition \lcx{\eqref{Cond.1}} can be written as
\[ \limit{n}(1-\rho_n) \ln n = \frac{2\lambda^2}{\tau}.\]
\bH{Note in passing that for the} Gaussian case, which corresponds
to $\tau =2$, the above asymptotic condition reduces to \bH{the
H\"usler-Reiss}  {condition \pE{\eqref{HRcond}}}.

 {\exxa (Weibull case)} Consider a triangular array as in \netheo{T2}
with almost surely positive random radius $R$ being \cL{b}eta
distributed
with parameters $a, b > 0$, thus %(cf. \cite{HashorvaLP13})
\[ \bar F(1- x) =  {\frac{\Gamma(a+b)}{b\Gamma(a)
\Gamma(b)}}x^{b}\left(1-\frac{b(a-1)x}{(b+1)}(1+o(1))\right),\quad
x\downarrow 0.
\] \rH{Consequently}, condition \eqref{Weibull} holds with $\alpha =
b$ and thus \lcx{\eqref{Cond.11} \pE{is satisfied if further}}
\[  \lim_{n\to\IF} (1-\rho_n) \fracl {b\Gamma(a)
\Gamma(b)}{n \Gamma(a+b)} ^ {-1/(1/2 + b )} = 2\lambda^2 \]
\pE{is valid.}
% \[ 1-\rho_n \sim 2\lambda^2 \fracl {b B(a, b)}{n} ^ {1/(1/2 + b )}, \quad \rH{n\to \IF}. \]

\section{Proofs}\label{sec4}
\lcx{We present first a lemma.} { Its proof is deferred to Appendix.} \\
For any $x\neq0, y\in\R$ and $\rho \in (0,1]$ denote
\BQN\label{Beta} \beta = \beta(x, y, \rho) = \arctan \fracl{y/x -
\rho}{\sqrt{1-\rho^2}},  \quad\psi = \arccos \rho \in[0, \pi/2).\EQN
\hH{In the following $c_n \sim d_n$ means $\limit{n} c_n/d_n=1$ for $c_n, d_n, n\ge 1$ given positive constants, and
instead of $B(y;a,b)$ we shall write {$B(y)$ simply} for the beta distribution with parameters $1/2,1/2$.} % simply $B(y)$. }\\

 \BL \label{L1} Let $(X, Y) \equaldis R (\cos \Theta, \sin \Theta)$ with independent random variables $R= \sqrt{X^2 +
 Y^2}
> 0$ almost surely and $\Theta \in (-\pi, \pi)$. Then for $x, y>0$ and $\psi, \beta$ given by
\eqref{Beta}{, we have}
% $a)$ For $\rho \in (0,1]$
\[ \Pk{X
> x, \rho |X| + \sqrt{1 - \rho^2} Y
> y}  =  \Pk{ R > \frac x {\cos\Theta}, \Theta \in( \beta , \frac \pi 2) } +
\Pk{ R > \frac y {\cos(\Theta - \psi)}, \Theta  \in ( -\frac \pi 2 +
\psi, \beta) } \]
 and
\BQNY  \Pk{ \rho |X| + \sqrt{1 - \rho^2} Y
> x } & = & \Pk{R > \frac x{\cos(\Theta - \psi)}, \Theta - \psi \in
(-\frac\pi 2, \frac \pi 2 -\psi) } \\
& \quad & + \Pk{ R > -\frac{x}{\cos(\Theta + \psi)}, \Theta \in
(\frac\pi 2 , \pi ) \cup (- \pi, -\frac\pi 2 - \psi )}. \EQNY
% $b)$ Let $\rho = \rho_1 \rho_2 + \sqrt{1-\rho_1^2
% }\sqrt{1-\rho_2^2}$ for given $\rho_1, \rho_2 \in (0, 1]$ and
% $\rho_1 \le \rho_2$, denote $\psi_j = \arccos \rho_j, j = 1, 2.$
% \begin{align*}& \quad  \Pk{\rho_1
% |X| + \sqrt{1 - \rho_1^2} Y > x, \rho_2 |X| + \sqrt{1 - \rho_2^2} Y
% >y} \\ & =  \Pk{ R > \frac x{\cos (\Theta - \psi_1)}, \Theta - \psi_1 \in (
% - \frac \pi 2,  \min(\frac\pi 2 - \psi_1, -\beta) }  + \Pk{ R >
% \frac y{ \cos (\Theta - \psi_2)}, \Theta - \psi_2 \in ( \psi
% - \beta,  \frac\pi 2 - \psi_2) }\\
% & \quad  + \Pk{ R > \frac{x}{ -\cos(\Theta + \psi_1)}, \Theta +
% \psi_1 \in (\max(\psi_1, \beta)-\pi , -\frac\pi 2 )\cup (\max(\beta
% + \pi, \frac\pi 2 +\psi_1), \pi + \psi_1)}
% \\  &\quad + \Pk{ R
% > \frac{y}{ -\cos(\Theta + \psi_2)}, \Theta + \psi_2 \in (\psi_2 - \pi,  \beta - \psi - \pi )
% \cup (\frac\pi 2 +\psi_2, \pi + \min( \psi_2, \beta - \psi))}.
% \end{align*}
 \EL
%\hh{The proof of \nelem{L1} is postponed to Appendix.}
 \prooftheo{T1}   By \nelem{L1} putting $\psi_n = \arccos \rho_n$ and $ u_n(y) = a_n y + b_n$ we obtain
\BQNY  \lefteqn{n\Pk{ \rho_n |X|+ \sqrt{1-\rho_n^2} Y > u_n(y) } }
%=
%n\Pk{ R (\rho_n |\cos \Theta| +
%\sqrt{1 - \rho_n^2} \sin \Theta) > u_n(y) } \\
\\&=& \frac{n}{2\pi} \left( \int_{-\pi/2}^0 \Pk{R \cos \theta
> u_n(y)} \, d \theta + \int_0^{\pi/2 -\psi_n }\Pk{ R \cos\theta
> u_n(y) } \,d\theta \right.\\
&\quad& \quad \left. + \int_{\psi_n + \pi/2}^{\pi}\Pk{ R
|\cos\theta|
> u_n(y) } \,d\theta  +\int_{\psi_n - \pi}^{-\pi/2} \Pk{ R
|\cos\theta|
> u_n(y) } \,d\theta + \int_{ \pi}^{\pi + \psi_n} \Pk{ R
|\cos\theta|
> u_n(y) } \,d\theta  \right) \\
& = & \frac{2n}{2\pi} \int_0^{\pi/2} \Pk{R \cos \theta
> u_n(y)} \, d \theta + \frac{2n}{2\pi}  \int_0^{\pi/2 -\psi_n }\Pk{ R \cos\theta
> u_n(y) } \,d\theta  =: A_n + B_n.\EQNY
It follows from
\eqref{F and G} that
\BQNY A_n = \frac{n}{2} \int_0^1 \Pk{R> \frac{u_n(y)}{\sqrt{1
-s}}}\, d B(s) = n \bar G(u_n(y)) \to  e^{-y}, \quad n\to \IF. \EQNY
Further, since
$\rho_n \le 1$ and $\rho_n
> \epsilon_0$ for sufficiently large $n$ and any given $\epsilon_0
\in (0,1)$,
\BQNY & & \int_0^{1} \Pk{R> \frac{u_n(y)}{\sqrt{1
-s}}}\, d B(s) \ge \int_0^{\rho_n^2} \Pk{R> \frac{u_n(y)}{\sqrt{1
-s}}}\, d
B(s)  = \pE{\frac{2 B_n}{n}}\\
&\quad &\quad  \ge
 \int_0^{\epsilon_0^2} \Pk{R>
\frac{u_n(y)}{\sqrt{1 -s}}}\, d B(s) \sim \int_0^{1} \Pk{R>
\frac{u_n(y)}{\sqrt{1 -s}}}\, d B(s), \quad n\to \IF,
\EQNY
\bH{where the last step above follows by} Proposition 12.2.1 %and Theorem 12.3.1
in \cite{Berman92}. Consequently,
\[ \limit{n}n \Pk{ \rho_n |X|+ \sqrt{1-\rho_n^2} Y > u_n(y)} = 2 e^{-y}.   \]
and thus \eqref{T1r1} follows. %Therefore, we complete the proof of $a)$.

Next we assume first that $x_F = \IF$, which implies that both
$u_n(x)$ and $u_n(y)$ tend to infinity as $n\to\IF$. By \nelem{L1}
\BQNY
n\Pk{X > u_n(x), \rho_n |X|+ \sqrt{1-\rho_n^2} Y >
u_n(y) } & = & \frac{n}{2\pi} \int_{\beta_n}^{\pi/2} \bar
F\fracl{u_n(x)}{\cos \theta}\, d\theta +
 \frac{n}{2\pi}
\int_{-\pi/2}^{\beta_n - \psi_n} \bar F\fracl{u_n(y)}{\cos \theta}\,
d\theta \\
&  =: & I_n + J_n,
\EQNY
where
\[ \beta_n = \arctan\left( \frac{u_n(y) / u_n(x) - \rho_n}{ \sqrt{1 - \rho^2_n}}\right), \quad \psi_n = \arccos \rho_n \in [0, \pi/2), \]
with
\BQNY
\tan (\beta_n) = \frac{u_n(y) / u_n(x) - \rho_n}{ \sqrt{1 -
\rho^2_n}} = \left( \frac{y-x + (1-\rho_n)x}{(1 - \rho_n) b_n w(
b_n)} + 1 \right) \sqrt{ \frac{1-\rho_n}{1+\rho_n}} \left( 1+
\frac{x}{b_n w( b_n)}\right)^{-1}.
 \EQNY
The asymptotic behaviours of integrals similar to $I_n$ and $J_n$ are derived in several contributions, see e.g.,
\cite{Berman92,Hashorva05S,Hashorva06,Hashorva08,Polar92,convex}.   The idea is to transform both the integrand and the distribution function (in our case $B(y)$) so that the integrand converges locally uniformly
and so does the sequence of distribution functions. Below we shall treat separately two
cases $(i)\ \lambda \in (0, \IF)$ \lcx{and} $(ii) \ \lambda = 0$. Denote $\nu_n = \nu(b_n) =
\sqrt{b_n w(b_n)}$.

$(i) $ If $\hH{\lambda_{x,y}}:=\lambda + (y - x)/(2\lambda) >0$,
then  $\beta_n > 0$ {holds} for sufficiently large  $n$. Thus, for
any given $\epsilon \in(0,1)$ and all $n$ \hH{large} \BQN\nonumber
I_n & = & \frac{n}{4} \int_{(\sin \beta_n)^2} ^1 \bar F
\fracl{u_n(x)}{\sqrt {1 - s }}\, d B(s)\\
\nonumber & \sim & \frac n{4\nu_n}\int_{(\sin\beta_n)^2
\nu^2_n}^{\epsilon^2 \nu^2_n} \bar F\fracl{\lcx{u_n(x)}}{\sqrt{1 -s/
\nu^2_n}}\, d \nu_n B(s/
\nu^2_n)\\
 \label{J1}& \sim & \left( 1 - \Phi \left(  \hH{\lambda_{x,y}}
\right) \right) e^{-x}
 \EQN
{s}ince $w( u_n(x) ) \sim w(b_n),  b_n w(b_n) \to \IF$ as $n\to
\IF$ and \BQNY ( \sin\beta_n )^2\nu_n^2  \sim  \beta^2_n \nu_n^2
\sim \left( \frac{y-x + (1-\rho_n)x}{(1 - \rho_n) b_n w( b_n)} + 1
\right)^2 \frac{1-\rho_n}{1+\rho_n} b_n w(b_n) \sim
\left(\hH{\lambda_{x,y}}\right)^2.
\EQNY
Similarly, \eqref{J1} holds
for $\hH{\lambda_{x,y}}\le 0$. Further, $\psi_n \sim \sqrt{2(1 -
\rho_n)} \to 0$ as $n\to\IF$ and
 \BQNY (\sin (\beta_n - \psi_n))^2 \nu_n^2 & \sim &
(\beta_n - \psi_n)^2 b_n w(b_n) \sim  b_n w(b_n) \left (
\frac{u_n(y) / u_n(x) - \rho_n}{\sqrt{1 - \rho_n^2}}-
\sqrt{2(1-\rho_n)}
\right)^2\\
& \sim & \left( \hH{\lambda_{y,x}}\right)^2.
\EQNY
Hence, it follows that for both two cases $\hH{\lambda_{x,y}}
> 2\lambda$ and $\hH{\lambda_{x,y}}
\le 2\lambda$ we have
 \BQN\label{J2}
 \hH{\limit{n}}J_n =\left( 1 - \Phi \left(
\hH{\lambda_{y,x}} \right) \right) e^{-y},
\EQN
% implying
{which implies that for  $(M_{n1}, M_{n2}):= (\max_{1\le j\le
n}X_{jn}, \max_{1\le j\le n}Z_{jn})$ and $(X_{jn}, Z_{jn}), j\le n$
given as in \eqref{H-skew}} \BQNY -\limit{n}n \ln \Pk{ M_{n1} \le
u_n(x), M_{n2} \le u_n(y)} = e^{-y} + \Phi \left(
\hH{\lambda_{x,y}}\right) e^{-x} + \Phi \left(
\hH{\lambda_{y,x}}\right)e^{-y}. \EQNY $(ii)$ For $\lambda = 0$,
note that as $n\to\IF$ \BQNY \label{Lim.1} (\sin\beta_n ) \nu_n & =
& \frac{(u_n(y) - \rho_n u_n(x)) \sqrt{ b_n
w(b_n)}}{\sqrt{(1-\rho_n^2) (u_n(x))^2 + (u_n(y) - \rho_n u_n(x))^2
}}\\
\nonumber & \sim &  \frac{(y- \rho_n x) \sqrt{ b_n w(b_n)}}{ \sqrt{2
(1-\rho_n)(b_n w(b_n) + x)^2 + (y- \rho_n x)^2}} \to \pm
\IF,\\
( \sin (\beta_n - \psi_n ) ) \nu_n  \label{Lim.2} & = & \frac{ (1 -
\rho_n) b_n w(b_n) + (y - \rho_n x) - (1 - \rho_n^2) (b_n w(b_n) +
x)}{\sqrt{(1 - \rho_n^2) (b_n w(b_n) + x)^2 + ((1 - \rho_n) b_n
w(b_n) + (y - \rho_n x))^2 }} \sqrt{b_n
w(b_n)} \\
\nonumber & \to & \pm \IF,
\EQNY
where the sign above depends on $y > x$ and $y
\le x$, respectively. Thus using
similar arguments {as above for (i)} it follows that both \eqref{J1} and \eqref{J2} hold also for $\lambda=0$. \\
If the upper endpoint $x_F \in (0, \IF)$, then $\bar F(x
/\cos\theta) = 0$ for $x/ \cos\theta > x_F$, thus one can substitute
the upper limits of the integrals above accordingly {and obtain
the results.} \hH{Hence} the proof is complete. \QED

\COM{ \prooftheo{T2} Without loss of generality, we assume that
$\rho_{1n} \le \rho_{2n}$, and thus denote
\BQNY \nonumber & &\psi_{jn} = \arccos \rho_{jn}, \quad \psi_n =
\psi_{1n} -
\psi_{2n}, \quad \rho_n = \cos \psi_n = \rho_{1n}\rho_{2n} + \sqrt{1 - \rho_{1n}^2}\sqrt{1 - \rho_{2n}^2}\\
  & & \beta_n  =  \arctan \fracl{u_n(y) / u_n(x) - \cos \psi_n }{\sin
\psi_n} = \arctan \fracl{u_n(y) / u_n(x) - \rho_n }{\sqrt {1 -
\rho^2_n}}.\EQNY

Note that $\psi_{jn} \sim \sqrt{ 2(1 - \rho_{jn})} \to 0$ as
$\rho_{jn} \to 1$, thus $\psi_n \to 0$ and further by \eqref{Cond.2}

\BQN \label{Diff.Lambda} (1- \rho_n )b_n w(b_n) \sim \frac{\psi_n^2
b_n w(b_n)}{2} = \frac{\Big( \sqrt{ 2(1 - \rho_{1n})} (1+o(1)) -
\sqrt{ 2(1 - \rho_{2n})} (1+o(1)) \Big)^2 b_n w(b_n)}{2} \sim
2\lambda^2 ,  \EQN

with $\lambda = |\lambda_1 - \lambda_2|.$

 Next, % denote $ Z_{jn} = \rho_{jn} |X| + \sqrt{1- \rho_{jn}^2} Y $  for $j =1, 2$.
by \nelem{L1}
\BQNY n \Pk{Z_{1n}
> u_n(x), Z_{2n} > u_n(y)} & = & \frac n {2\pi} \left( \int_{-\pi/2}^{ - \beta_n} \bar F\fracl{u_n(x)}{\cos\theta} \,
d\theta +\quad \int_{\psi_n- \beta_n}^{\pi/2 - \psi_{2n}}  \bar
F\fracl{u_n(y)}{\cos\theta} \, d\theta  \right.
 \\
&\quad & \quad + \int_{\max(\beta_n, \psi_{1n})-\pi}^{ - \pi/2} \bar F\fracl{u_n(x)}{-\cos\theta} \, d\theta + \int_{\beta_n + \pi}^{\max(\beta_n, \psi_{1n}) + \pi} \bar F\fracl{u_n(x)}{-\cos\theta} \, d\theta \\
& \quad & \quad \left. + \int_{\psi_{2n}- \pi}^{\max(\psi_{2n},
\beta_n - \psi_n) - \pi}  \bar F\fracl{u_n(y)}{-\cos\theta} \,
d\theta
+  \int_{\psi_{2n}+ \pi/2}^{\min(\psi_{2n}, \beta_n - \psi_n) + \pi}  \bar F\fracl{u_n(y)}{-\cos\theta} \, d\theta \right)\\
& = & \frac{n}{2\pi} \left( \int_{\beta_n}^{\pi/2} +
\int_{\max(\beta_n, \psi_{1n})}^{\pi/2} + \int_{\beta_n}^{\max(\beta_n, \psi_{1n})} \bar F\fracl{u_n(x)}{\cos\theta} \, d\theta \right. \\
& \quad & \quad + \left. \int_{\psi_n- \beta_n}^{\pi/2 - \psi_{2n}}
 + \int_{\psi_{2n}}^{\max(\psi_{2n}, \beta_n - \psi_n)} + \int_{-\min(\psi_{2n}, \beta_n -\psi_{n})}^{\pi/2 -\psi_{2n}} \bar F\fracl{u_n(y)}{\cos\theta} \, d\theta\right)\\
 & =: & \sum _{j=1}^6 L_j.\EQNY

Note that for $\lambda \in(0, \IF)$ \BQNY \beta_n \nu_n \to \lambda
+ \frac{ y - x}{2\lambda}, \quad \psi_{jn} \nu_n \to 2\lambda_j, j =
1, 2.\EQNY
One can show for both two cases $\lambda +
(y-x)/(2\lambda) > 2\lambda_1$ and $\lambda + (y-x)/(2\lambda) \le
2\lambda_1$ that (cf. \eqref{J1})

\[I_n + J_n + K_n \sim 2e^{-x} \left( 1 - \Phi \left(  \lambda + \frac{y-x}{2\lambda}
\right) \right). \]

Similar arguments of \eqref{J2} yield that

\[I_4 + I_5 + I_6 \sim 2e^{-y} \left( 1 - \Phi \left(  \lambda + \frac{x-y}{2\lambda}
\right) \right). \]

The other cases for $\lambda = 0$ or $x_F<\IF$ can be established
accordingly. Therefore, the desired results are obtained. \QED }

\prooftheo{T2} $a)$ First recall that
\BQNY
z_n = \arccos \rho_n \in
[0, \pi/2), \quad \psi_n = \arccos (1 - u_n) \sim \sqrt{2u_n}\to 0,
\quad  {n \to \IF}. \EQNY By \nelem{L1} and the fact that $x_F=1$ we
have \BQNY \Pk{\rho_n |X| + \sqrt{1 - \rho_n^2} Y > 1 - u_n}
%& = & \Pk{R \cos (\Theta - z_n) > 1 - u_n, \Theta - z_n \in
%(-\frac \pi 2, \frac \pi 2 - z_n)}
%\\
%& \quad & + \Pk{R \cos (\Theta + z_n) <- (1 - u_n), \Theta + z_n\in
%(\frac \pi
%2 + z_n, \pi + z_n ) \cup ( z_n-\pi, -\frac \pi 2 )} \\
& = & \Pk{R \cos (\Theta - z_n) > 1 - u_n, \Theta - z_n \in
(-\psi_n, \psi_n)}
\\
& \quad & + \Pk{R \cos (\Theta + z_n) <- (1 - u_n), \Theta + z_n \in
(\pi
 - \psi_n, \pi + \min (z_n, \psi_n) )} \\
 & \quad & + \Pk{R \cos (\Theta + z_n) <- (1 - u_n), \Theta + z_n \in (z_n - \pi, \psi_n - \pi)}\\
 & =: & I_n + J_n + K_n.
 %\\
\EQNY
%using the upper endpoint of $R$ is 1.
\aH{Next, for all large $n$}
 \BQNY
 I_n
& = & \frac{2}{2\pi} \int_0^{\psi_n} \bar F\fracl {1 - u_n} { \cos
\theta} \, d\theta = \frac{1}{2} \int_0^{(\sin\psi_n)^2} \bar F
\fracl{ 1 - u_n}{ \sqrt{1- s}} \,
d B(s)\\
& = &  \frac{1}{2} \int_0^{(\sin\psi_n)^2 / (2 u_n)} \bar F\fracl{ 1
- u_n} { \sqrt{1- 2 u_ns}} \, d B(2 u_ns) \\
& = & \frac{\sqrt{2 u_n }\bar F(1 - u_n) }{2}
\int_0^{(\sin\psi_n)^2/ (2 u_n)} \frac{ \bar F( (1 - u_n )/ \sqrt{1-
2 u_ns}) }{\bar F(1 - u_n)} \, d \frac{B(2 u_ns)}{\sqrt{2 u_n }}.
\EQNY
By Theorem 12.3.3 in \cite{Berman92} (see also Theorem 2.1 in \cite{Hashorva08})
\[ I_n \sim \mathcal {I}_\alpha u_n ^{1/2} \bar F(1 - u_n),  \quad n\to \IF,\]
with $\mathcal {I}_\alpha$ {given by} \eqref{I(alpha)}.  For
$J_n$ and $K_n$, we consider first that $2\lambda^2
> 1$, then  $z_n > \psi_n$ for sufficiently large $n$. Thus
$K_n = 0$ and \BQNY J_n & = & \Pk{R \cos (\Theta + z_n) <- (1 -
u_n), \Theta + z_n \in (\pi
 - \psi_n, \pi +  \psi_n )} \\
& = & \Pk{R \cos (\Theta + z_n - \pi) >  1 - u_n, \Theta + z_n -\pi
\in (
 - \psi_n,  \psi_n )}  = I_n.
 \EQNY
Similarly for $2\lambda^2 < 1$ and $2\lambda^2 = 1$ \hH{it follows that} $J_n + K_n = I_n$. %Thus \netheo{T2} holds with $x = -1$.
Hence the proof of $a)$ follows by utilising further \eqref{Weibull}.

$b)$ For any $x,y$ negative the case proved in $a)$ above implies
{that}
 \BQNY
\limit{n} n\Pk{ \rho_n |X| + \sqrt{1 - \rho_n^2} Y > 1 + u_n  y} =
2\mathcal {I}_\alpha  \abs y^{1/2 +
 \alpha}
 \EQNY
 and
\BQNY \limit{n}n\Pk{ X > 1 + u_n  x } = \pE{\limit{n}}\frac n 2 \Pk{|X| > 1 + u_n
x} = \mathcal {I}_\alpha  \abs x^{1/2 +
 \alpha}. \EQNY % So
{It thus} remains to deal with
\[ A_n := \Pk{ X> 1 + u_n  x, \rho_n |X| + \sqrt{1 - \rho_n^2} Y > 1 +u_n  y}. \]
By \nelem{L1} \pE{we have}
\begin{align*} A_n & =   \Pk{ R \cos\Theta  > 1 + u_n x, \Theta \in(
\beta_n , \pi/2) } + \Pk{ R \cos(\Theta - z_n) > 1 + u_n y, \Theta -
z_n
\in ( -\pi/2, \beta_n - z_n) }  \\
& =  \Pk{ R > \frac{ 1 + u_n x } {\cos\Theta} , \Theta \in( \max(
-\psi_{1n}, \beta_n), \psi_{1n} )} + \Pk{ R > \frac{1 + u_n
y}{\cos(\Theta - z_n)}, \Theta - z_n \in (
-\psi_{2n}, \min(\beta_n - z_n, \psi_{2n}) ) }\\
& =:  A_{1n} + A_{2n},
 \end{align*}
 where $z_n = \arccos \rho_n \sim 2\lambda \sqrt{ u_n}$ as $n\to\IF$ and
\BQNY & \quad & \psi_{1n} = \arccos (1 + u_n x) \sim \sqrt {2u_n
\abs x},
\quad \psi_{2n} = \arccos (1 + u_n y) \sim \sqrt{2u_n \abs y},\\
& \quad & \beta_n = \arctan\fracl{(1+ u_n y)/ (1+ u_n x) -
\rho_n}{\sqrt{1 - \rho_n^2}} \sim \left(\lambda + \frac{y-
x}{2\lambda} \right)\sqrt{ u_n}. \EQNY If $0 <
\hH{\lambda_{x,y}:=}\lambda + \frac{y-x}{2\lambda} < \sqrt{2\abs
x}$, then $0 < \beta_n < \psi_{1n}$ holds for sufficiently large
$n${. Hence}%, hence
\BQN \nonumber n A_{1n} & = & \frac n{2\pi}
\int_{\beta_n}^{\psi_{1n}} \bar F\fracl{1+u_n x}{\cos\theta}\,
d\theta = \frac n 4 \int_{(\sin \beta_n)^2}^{(\sin \psi_{1n} )^2}
\bar F\fracl{1 + u_n x}{\sqrt {1 - s}}\, d B(s) \\
\nonumber & = & \frac n 4 \int_{(\sin \beta_n)^2 / (2u_n \abs
x)}^{(\sin \psi_{1n} )^2 / (2u_n \abs x)}
\bar F\fracl{1 + u_n x}{\sqrt {1 + 2u_n x s}}\, d B(2 u_n \abs x s)\\
\nonumber & \sim & \abs x^{1/2 + \alpha}(nu_n^{1/2} \bar F(1-u_n))
 \frac{\sqrt 2}{2\pi}\int_{\hH{\lambda_{x,y}} /
\sqrt{2\abs x}}^1
(1-s^2)^\alpha\, ds \\
\label{Asym.4} & \sim & \abs x^{1/2 + \alpha}\frac{\sqrt
2}{2\pi}\int_{\hH{\lambda_{x,y}}/ \sqrt{2\abs x}}^1
(1-s^2)^\alpha\, ds = \abs x^{1/2 + \alpha} \mathcal {I}_\alpha
\left( 1 - \Upsilon_\alpha \fracl{\hH{\lambda_{x,y}}}{\sqrt{2\abs x}} \right),
 \EQN
where $\mathcal {I}_\alpha$ and $\Upsilon_\alpha (\cdot)$ are given
by \eqref{I(alpha)}.

 {If $- \sqrt{2\abs x}<\hH{\lambda_{x,y}}< 0$,} then
$-\psi_{1n} < \beta_n \le 0$ holds for all large $n$% implying
{. Hence}
 \BQNY
\nonumber n A_{1n} & = & \frac n{2\pi} \left(\int_0^{-\beta_n}\bar
F\fracl{1+u_n x}{\cos\theta}\, d\theta + \int_0^{\psi_{1n}} \bar
F\fracl{1+u_n x}{\cos\theta}\,
d\theta\right) \\
& \sim & \abs x^{1/2 + \alpha}(nu_n^{1/2} \bar F(1-u_n))
 \frac{\sqrt 2}{2\pi}\left(\int_0^{-\hH{\lambda_{x,y}}/
\sqrt{2\abs x}}(1-s^2)^\alpha\, ds + \int_0^1
(1-s^2)^\alpha\, ds \right) \\
 & \sim & \abs x^{1/2 + \alpha}\frac{\sqrt
2}{2\pi}\int_{\hH{\lambda_{x,y}} / \sqrt{2\abs x}}^1
(1-s^2)^\alpha\, ds = \abs x^{1/2 + \alpha} \mathcal {I}_\alpha
\left( 1 - \Upsilon_\alpha \fracl{\hH{\lambda_{x,y}}}{\sqrt{2\abs x}} \right).
 \EQNY
Similarly, \eqref{Asym.4} holds for the other three cases: $\hH{\lambda_{x,y}}= - \sqrt{2\abs x},$ $\hH{\lambda_{x,y}}= 0$ and $\hH{\lambda_{x,y}}\ge
\sqrt{2\abs x}$, respectively.

For $A_{2n}$, if $0 < \hH{\lambda_{y,x}} \hh{:=}%=:
\lambda + \frac{x-y}{2\lambda} < \sqrt{2\abs y}$, then $ -
\psi_{2n}< \beta_n - z_n < 0$  {holds} for sufficiently large $n$% and therefore
{. Hence} \BQNY n A_{2n} & = & \frac n{2\pi}
\int_{-\psi_{2n}}^{\beta_n - z_n} \bar F\fracl{1+u_n
y}{\cos\theta}\, d\theta = \frac n 4 \int_{(\sin (z_n
-\beta_n))^2}^{(\sin \psi_{2n} )^2}
\bar F\fracl{1 + u_n y}{\sqrt {1 - s}}\, d B(s) \\
& = & \frac n 4 \int_{(\sin (z_n -\beta_n))^2 / (2u_n \abs
y)}^{(\sin \psi_{2n} )^2 / (2u_n \abs y)}
\bar F\fracl{1 + u_n y}{\sqrt {1 + 2u_n y s}}\, d B(2 u_n \abs y s)\\
& \sim & \abs y^{1/2 + \alpha}(nu_n^{1/2} \bar F(1-u_n))
\frac{\sqrt 2}{2\pi}\int_{\hH{\lambda_{y,x}} /
\sqrt{2\abs y}}^1
(1-s^2)^\alpha\, ds \\
& \sim & \abs y^{1/2 + \alpha} \frac{\sqrt 2}{2\pi}\int_{\hH{\lambda_{y,x}} / \sqrt{2\abs y}}^1 (1-s^2)^\alpha\, ds = \abs
y^{1/2 + \alpha}\mathcal {I}_\alpha  \left( 1 - \Upsilon_\alpha
\fracl{\hH{\lambda_{y,x}}}{\sqrt{2\abs y}} \right).
 \EQNY
 {The other three cases $\hH{\lambda_{y,x}}\le - \sqrt{2\abs y},  - \sqrt{2\abs y} <\hH{\lambda_{y,x}}\le 0$ and $\hH{\lambda_{y,x}}\ge\sqrt{2\abs y}$} follow with similar arguments as above establishing thus the
claim in $b)${. Consequently,} the proof is complete. \QED

\COM{ \prooftheo{T4} Without loss of generality, we assume that
$\rho_{1n} \le \rho_{2n}$, thus $\lambda := \lambda_1 - \lambda_2
\ge 0$. Denote $\rho_n = \rho_{1n}\rho_{2n} + \sqrt{1 - \rho_{1n}^2}
\sqrt{1 - \rho_{2n}^2}$ and thus it follows from \eqref{Cond.2} and
similar arguments of \eqref{Diff.Lambda} that
\[\frac{1 - \rho_n}{u_n} \to  2\lambda^2.\]
Further \BQNY & \quad & z_{jn} := \arccos \rho_{jn}\sim 2\lambda_j
\sqrt{u_n}, \quad \beta_n := \arctan \fracl{(1+u_ny)/(1+u_nx) -
\rho_n}{\sqrt{1-\rho_n^2}} \sim \left(\lambda + \frac{x -
y}{2\lambda} \right) \sqrt{u_n}, \\
&\quad & z_n := z_{1n} - z_{2n}\sim 2\lambda \sqrt{u_n}, \quad
\psi_{1n} := \arccos (1+ u_n x) \sim \sqrt{2r_nx}, \quad \psi_{2n}
:= \arccos (1+ u_n y)\sim \sqrt{2r_ny}.\EQNY
Recall that the upper
endpoint of $R$ is 1, by \nelem{L1}
\begin{align*} &  \Pk{\rho_{1n} |X| + \sqrt{1 - \rho_{1n}^2} Y > 1 + u_n x, \rho_{2n} |X| + \sqrt{1 - \rho_{2n}^2} Y > 1 + u_n y}\\
& =  \Pk{ R  > \frac{1+ u_n x}{\cos (\Theta - z_{1n})}, \Theta -
z_{1n} \in ( - \psi_{1n},  \min(\psi_{1n}, -\beta_n) } + \Pk{ R  >
\frac{1+u_ny}{\cos (\Theta - z_{2n})}, \Theta - z_{2n} \in ( \max(
-\psi_{2n}, z_n
- \beta_n ),  \psi_{2n}) }\\
& \quad + \Pk{ R > \frac{-(1+u_nx)}{ \cos(\Theta + z_{1n})}, \Theta
+ z_{1n} \in (\max(-\psi_{1n}, z_{1n}, \beta_n)-\pi , \psi_{1n} -
\pi)\cup (\max(-\psi_{1n}, \beta_n)+ \pi, \pi +
\min(\psi_{1n},z_{1n}))}
 \\
 & \quad + \Pk{ R
> \frac{-(1+u_n y)}{ \cos(\Theta + z_{2n})}, \Theta + z_{2n} \in (\max(- \psi_{2n}, z_{2n}) - \pi,  \min(\psi_{2n}, \beta_n - z_n) - \pi )
\cup (\pi - \psi_{2n}, \pi + \min( z_{2n}, \psi_{2n}, \beta_n -
z_n))}\\
& =: \sum_{j = 1}^4 I_j.\end{align*} Next we treat the four terms in
turn.

For $I_n$, if $ \sqrt{2x} < - (\lambda + (x-y)/(2\lambda))$, then
$\psi_{1n} < - \beta_n$ for sufficiently large $n$, and thus
\BQNY n
I_n = \frac{2n}{2\pi} \int_0^{\psi_{1n}} \Pk{R > \frac{1 +u_n
x}{\cos\theta}} \, d\theta \sim x^{1/2 + \alpha} \mathcal \mathcal
{I}_\alpha \COM{= x^{1/2 + \alpha} \mathcal {I}_\alpha \left( 1 -
\Upsilon_\alpha \fracl{\lambda +
(y-x)/(2\lambda)}{\sqrt{2x}}\right)}.\EQNY

If $ \sqrt{2x} > - (\lambda + (x-y)/(2\lambda))$, then $\psi_{1n} >
- \beta_n \ge 0$ or $\psi_{1n} >0 \lcx{>} - \beta_n$, one can show
(cf. \eqref{Asym.4})

\[n I_n \sim x^{1/2 + \alpha} \mathcal {I}_\alpha  \left( 1 -
\Upsilon_\alpha \fracl{\lambda +
(x-y)/(2\lambda)}{\sqrt{2x}}\right).\]

Similarly, one can show the same results hold for $ \sqrt{2x} = -
(\lambda + (x-y)/(2\lambda))$.

For $K_n$, note that

\BQNY K_n = \frac{n}{2\pi} \left( \int_{\theta \in (\max(z_{1n},
\beta_n), \psi_{1n})} \Pk{R > \frac{1+u_nx}{\cos\theta} }\, d\theta
+ \int_{\theta\in (\max(-\psi_{1n}, \beta_n), \min(z_{1n},
\psi_{1n})) }\Pk{R
> \frac{1+u_nx}{\cos\theta} }\, d\theta \right). \EQNY

If $ \sqrt{2x} < \lambda + (x-y)/(2\lambda)$, then $\psi_{1n} <
\beta_n$ for sufficiently large $n$, and thus $K_n = 0$, and if $
\sqrt{2x}
> \lambda + (x-y)/(2\lambda)$, then $\psi_{1n} > \beta_n$  for sufficiently large $n$, one can
follow the arguments of $J_n$ for the three cases: $a)\ z_{1n}
\lcx{\le} \beta_n; \ b) \ z_{1n} \lcx{\ge} \psi_{1n}$ and $c)  \
\psi_{1n}
> z_{1n}
> \beta_n$ and obtain $K_n \sim I_n$.

Next we deal with $J_n$ and $I_4$, respectively. If $- \sqrt{2y} <
2\lambda - (\lambda + (y-x)/(2\lambda))$,  then $- \psi_{2n} < z_n -
\beta_n$ for sufficiently large $n$. Further one can follow the
similar arguments of $I_n$ according to $z_n > \beta_n$ and $z_n
\lcx{\le} \beta_n$ and obtain that

\BQNY J_n \sim y^{1/2 + \alpha} \mathcal {I}_\alpha \left( 1 -
\Upsilon_\alpha \fracl{\lambda +
(y-x)/(2\lambda)}{\sqrt{2y}}\right).\EQNY

Similar argument for $I_4$ yields that $I_4\sim J_n$. Consequently,

\BQNY & \quad & \Pk{\rho_{1n} |X| + \sqrt{1 - \rho_{1n}^2} Y > 1 +
u_n x,
\rho_{2n} |X| + \sqrt{1 - \rho_{2n}^2} Y > 1 + u_n y}\\
 & \quad & \sim 2 \abs x^{1/2 + \alpha} \mathcal {I}_\alpha  \left( 1 -
\Upsilon_\alpha \fracl{\lambda + (y-x)/(2\lambda)}{\sqrt{2\abs
x}}\right) + 2 \abs y^{1/2 + \alpha} \mathcal {I}_\alpha  \left( 1 -
\Upsilon_\alpha \fracl{\lambda + (x-y)/(2\lambda)}{\sqrt{2\abs
y}}\right). \EQNY

Combining with \netheo{T3} $a)$, the desired results are obtained.
\QED}

\textbf{Acknowledgments}. \hH{We would like to thank the referees for several comments and suggestions. C. Ling kindly acknowledges support from
the Swiss National Science Foundation Grant 200021-140633/1. Both authors
have been partially supported by the project RARE -318984 (an FP7  Marie Curie
IRSES Fellowship).}

\section*{Appendix}
\cL{The Appendix contains the proof of \nelem{L1} and the direct
verification of \eqref{eqAB} for the half-skew Gaussian case.}

 \prooflem{L1} Recall that  $\psi = \arccos \rho \in[0, \pi/2)$
and the representation given by \eqref{Def.1} implies \BQNY (X, \rho
|X| + \sqrt{1 - \rho^2} Y) \equaldis R(\cos\Theta, |\cos\Theta|\cos
\psi + \sin\Theta \sin \psi), \EQNY where $\Theta
%\sim U
\cL{\in}(-\pi, \pi)$, independent of the random radius $R$. Hence,
\BQNY \Theta_\rho  :=  \rho |\cos \Theta| + \sqrt{1 - \rho^2} \sin
\Theta = \left\{ \begin{array}{rl}
\cos (\Theta - \psi), & \Theta \in (-\pi/2, \pi/2), \\
          - \cos (\Theta + \psi), & \Theta\in (-\pi, - \pi/2) \cup (\pi/2,
          \pi).
\end{array}
\right. \EQNY Further, $\Theta_\rho > 0$ if and only if \BQNY \Theta_\rho =
\left\{
\begin{array}{rl}
          \cos (\Theta - \psi), & \Theta - \psi \in (-\pi/2, \pi/2 - \psi), \\
          - \cos (\Theta + \psi), & \Theta + \psi \in (\psi -\pi, - \pi/2) \cup (\psi + \pi/2,
            \psi + \pi).
        \end{array}
\right. \EQNY Therefore, for $x, y > 0$ \BQNY  \Pk{ \rho |X| +
\sqrt{1 - \rho^2} Y
> x } & = & \Pk{R > \frac x{\cos(\Theta - \psi)}, \Theta - \psi \in
(-\frac\pi 2, \frac \pi 2 -\psi) } \\
& \quad & + \Pk{ R > -\frac{x}{\cos(\Theta + \psi)}, \Theta \in
(\frac\pi 2 , \pi ) \cup (- \pi, -\frac\pi 2 - \psi )} \EQNY and
 \BQNY  \Pk{X > x, \rho |X| + \sqrt{1 - \rho^2} Y
>y} & = & \Pk{R\cos \Theta > x, R \cos( \Theta - \psi ) > y, \cos \Theta
> 0 } \\
& = & \Pk{R > \max\left( \frac x {\cos \Theta}, \frac y {\cos(
\Theta - \psi )} \right), \Theta \in( \psi - \frac\pi 2, \frac\pi 2) } \\
& = & \Pk{ R > \frac x {\cos\Theta}, \Theta \in( \beta , \frac\pi 2)
} + \Pk{R > \frac y {\cos(\Theta - \psi)}, \Theta  \in ( -\frac\pi 2
+ \psi, \beta) },
 \EQNY
where $\beta$ is the solution of $\cos( \theta - \psi ) / \cos\theta
= y/x$ with respect to $\theta \in (-\pi/2, \pi/2)$, i.e.,
\[ \beta := \arctan \fracl {y/x - \rho} {\lcx{\sqrt{1-\rho^2}}}.\]
\QED
 \COM{ $b)$ Similar to $a)$, we can proceed for $\rho := \rho_1
\rho_2 + \sqrt{1- \rho_1^2} \sqrt{1- \rho_2^2} \in (0, 1]$ and $\psi
= \arccos \rho = \psi_1 - \psi_2$ with $\psi_j = \arccos \rho_j, j=
1,2$
 \begin{align*} & \quad  \Pk{\rho_1 |X| + \sqrt{1 - \rho_1^2} Y
> x, \rho_2 |X| + \sqrt{1 - \rho_2^2} Y
>y} \\
& =  \Pk{ R \cos (\Theta - \psi_1) > x, R \cos (\Theta - \psi_2) >
y, \Theta \in (-\frac\pi 2 + \psi_1, \frac\pi 2)} \\
&\quad   + \Pk{ R \cos (\Theta + \psi_1) < -x, R \cos (\Theta +
\psi_2) < -y, \Theta \in (-\pi, -\frac\pi 2 - \psi_1) \cup (
\frac\pi 2, \pi) }
\\
& = \Pk{ R \cos (\Theta - \psi_1) > x, \Theta - \psi_1 \in ( - \frac
\pi 2,  - \beta)\cap (-\frac\pi 2, \frac\pi 2 - \psi_1) }   +\Pk{ R
\cos (\Theta - \psi_2) > y, \Theta - \psi_1 \in (- \beta, \frac \pi
2)\cap(-
\frac \pi 2,  \frac \pi 2 - \psi_1)} \\
&\quad  +  \Pk{ R \cos (\Theta + \psi_1) < -x,   \Theta + \psi_1
\in\Big( ( \psi_1 -\pi, -\frac\pi 2) \cup (  \psi_1 + \frac\pi 2, \psi_1 +\pi)\Big) \cap \Big(( \beta-\pi, -\frac\pi 2) \cup ( \beta + \pi,  \frac{3\pi} 2)\Big)}\\
&\quad  + \Pk{ R \cos (\Theta + \psi_2) < -y,   \Theta + \psi_2
\in\Big( ( \psi_2 -\pi, -\frac\pi 2) \cup (  \psi_2 + \frac\pi 2,
\psi_2 +\pi)\Big) \cap \Big(( \frac\pi 2 - \psi, \beta- \psi +\pi)
\cup ( -\frac{3\pi} 2- \psi, \beta - \psi- \pi )\Big)}.
  \end{align*}
Here $\beta$ is the solution of the equation $ \cos (\theta' -
\psi)/ \cos \theta'  = y / x$ with respect to $\theta' = \psi_1 -
\theta$, i.e.,
\[ \beta := \arctan \fracl {y/x - \rho} {\sqrt{1-\rho^2}}.\]
Consequently, the desired results are obtained.}

\cc{{\bf Proof of the claim in \eqref{eqAB}}}.

Recall that $X, Y$ are independent standard normal distributed and
$b_n = G^\leftarrow (1-1/n) \sim \sqrt{2\ln n}, a_n = 1/b_n$, thus
$\varphi(b_n) / b_n \sim 1/n$ with $\varphi(\cdot)$ the
density of an $N(0,1)$ df. Further (set
$u_n(y) = y/b_n + b_n$)
 \BQNY & \quad & n\Pk{ \rho_n |X|+ \sqrt{1-\rho_n^2} Y > u_n( y  ) } = 2n \int _0^\IF
\Pk{Y > \frac{ u_n(
 y  ) - \rho_n x}{\sqrt{1 - \rho_n^2}}} \, d
\Phi(x)\\
&=& 2 \int _{- b_n^2}^\IF  \frac{n\varphi(u_n(x))}{b_n}\Pk{Y > \frac{
u_n(  y  ) - \rho_n u_n(x) }{ \sqrt{1 - \rho_n^2}}} \, d x \\
& \sim & 2 \int _{- b_n^2}^\IF  \exp\left( -\left( x +
\frac{x^2}{2b_n^2}\right)\right) \left(  1-\Phi\left( \frac{ y
-x}{\sqrt{1-\rho_n^2} b_n} + \frac{(1-\rho_n)x}{\sqrt{1-\rho_n^2}
b_n} + \frac{1-\rho_n}{\sqrt{1-\rho_n^2}} b_n \right) \right) \, d x
\\
&=& 2  \int _0^{ b_n^2}   \exp\left(x - \frac{x^2}{2b_n^2}\right)
\left(  1-\Phi\left( \frac{ y +x}{\sqrt{1-\rho_n^2} b_n} +
\frac{(1-\rho_n)(b_n^2 - x)}{\sqrt{1-\rho_n^2} b_n} \right) \right)
\, d x
\\
&\quad & +2\int _0^{\IF}\exp\left( -\left( x +
\frac{x^2}{2b_n^2}\right)\right) \left(  1-\Phi\left( \frac{ y
-x}{\sqrt{1-\rho_n^2} b_n} + \frac{(1-\rho_n)x}{\sqrt{1-\rho_n^2}
b_n} + \frac{1-\rho_n}{\sqrt{1-\rho_n^2}} b_n \right) \right)\, d x
\\
& =: & 2(I_n + J_n). \EQNY

Note that
\BQNY I_n \le  \int _0^{ b_n^2}   e^x \left(  1-\Phi\left( \frac{ y
+x}{\sqrt{1-\rho_n^2} b_n}  \right) \right) \, d x \le \int _0^{
\IF} e^x \left(  1-\Phi\left( \frac{ y +x}{\sqrt{1-\rho_n^2} b_n}
\right) \right) \, d x < \IF, \EQNY
where the last inequality holds by the dominated convergence theorem. In
fact $\sqrt{1-\rho_n^2} b_n < 2\lambda + 1, y+ x > 2\lambda + 1$
hold for sufficiently large $n$ and $x$. \pE{Further since
$1-\Phi(x) < \varphi(x) / x, x>0$ we obtain}
\[  e^x \left(  1-\Phi\left( \frac{ y
+x}{\sqrt{1-\rho_n^2} b_n} \right) \right) < \frac{\varphi((y+x)/
(2\lambda + 1 ))}{(y+x)/ (2\lambda + 1 )}e^x  \le
\varphi\left(\frac{y+x}{2\lambda + 1 }\right) e^x.
\]
Consequently, as $n\to \IF$
 \BQNY I_n \sim
2 \int _0^\IF \left( 1-\Phi\left( \frac{ y + x}{2\lambda} + \lambda
\right) \right)e^{x} \, d x =: 2I, \quad J_n \sim 2 \int_0^\IF
\left( 1-\Phi\left( \frac{ y -x}{2\lambda} + \lambda \right)
\right)e^{-x} \, d x =: 2J.
 \EQNY
Partial integration for $I$ and $J$ yields that (cf.
\cite{HuslerR89})
\BQNY && I = -\left(  1-\Phi\left( \frac{ y
}{2\lambda} + \lambda \right) \right) + \frac{e^{ -y }}{2\lambda}
\int_0^\IF \varphi\left( \frac{ y  + x}{2\lambda} - \lambda
\right)\, d x = -\left( 1-\Phi\left( \frac{ y }{2\lambda} + \lambda
\right) \right) + e^{ -y
} \Phi\left( \lambda - \frac{ y }{2\lambda}\right)\\
 && J = 1 - \int_0^\IF \Phi\left( \frac{ y -x}{2\lambda} + \lambda
\right) e^{-x}  \, d x = 1 - \Phi\left( \frac{ y }{2\lambda} +
\lambda \right) + e^{ -y }\left( 1 -  \Phi\left( \lambda - \frac{ y
}{2\lambda}\right) \right).
\EQNY
Hence,
 \BQNY \limit{n}  n\Pk{ \rho_n |X|+ \sqrt{1-\rho_n^2} Y > u_n( y  ) } = 2 e^{-y}, \EQNY
i.e.,  \BQNY \limit{n} n\Pk{ \rho_n |X|+ \sqrt{1-\rho_n^2} Y > a_n y
+ (b_n + a_n \ln2) }
 =
e^{-y}. \EQNY
 Next, note that $ u_n(x)
>0$ for all $x\in \R$ and sufficiently large $n$
\BQN \nonumber n \Pk{X > u_n(x), \rho_n |X|+ \sqrt{1-\rho_n^2} Y >
u_n(y) } & = & n \int_{u_n(x)}^\IF \Pk{Y > \frac{u_n(y) - \rho_n
t}{\sqrt{1-\rho_n^2}}} \, d \Phi(t)
\\
\nonumber & = & \int_x^\IF \frac{n \varphi( u_n(t) )}{b_n} \left( 1-
\Phi\left( \frac{u_n(y) - \rho_n u_n(t)}{\sqrt{1-\rho_n^2}}
\right)\right) \,
dt \\
\nonumber & \lcx{\sim} & \int_x^\IF \frac{\varphi( u_n(t) )}{\varphi( b_n
)} \left( 1- \Phi\left( \frac{(y-t) +
(1-\rho_n)t}{\sqrt{1-\rho_n^2}b_n}+ \sqrt{
\frac{1 - \rho_n }{1 + \rho_n}} b_n \right)\right) \, dt \\
\label{Gauss} & \sim & \int_x^\IF  \left( 1- \Phi\left(
\frac{y-t}{2\lambda}+ \lambda \right)\right) e^{-t} \, dt
 \EQN
 by the dominated convergence theorem. It follows from
 \cite{HuslerR89} that the left-hand side of \eqref{Gauss} is
 asymptotically equal to
\[  e^{-x} + e^{-y} - \left( \Phi
\left( \lambda +\frac{x-y}{2\lambda}\right) e^{-y} + \Phi \left(
\lambda +\frac{y-x}{2\lambda}\right) e^{-x} \right),
\]
establishing the claim. \QED

\bibliographystyle{plain}
 \bibliography{extremvalueseAO}
\end{document}